\documentclass[11pt,a4paper]{article}
\usepackage{a4,amssymb,latexsym,amsmath}
\begin{document}


\newcounter{pkt}
\newenvironment{enum}{\setcounter{pkt}{0}
\begin{list}{\rm\alph{pkt})}{\usecounter{pkt}
\setlength{\topsep}{1ex}\setlength{\labelwidth}{0.5cm}
\setlength{\leftmargin}{1cm}\setlength{\labelsep}{0.25cm}
\setlength{\parsep}{-3pt}}}{\end{list}~\\[-6ex]}

\newcounter{punkt1}
\newenvironment{enum1}{\setcounter{punkt1}{0}
\begin{list}{\arabic{punkt1})}{\usecounter{punkt1}
\setlength{\topsep}{1ex}\setlength{\labelwidth}{0.6cm}
\setlength{\leftmargin}{1cm}\setlength{\labelsep}{0.25cm}
\setlength{\parsep}{-3pt}}}{\end{list}~\\[-6ex]}

\newcounter{punkt3}
\newenvironment{enum2}{\setcounter{punkt3}{0}
\begin{list}{(\roman{punkt3})}{\usecounter{punkt3}
\setlength{\topsep}{1ex}\setlength{\labelwidth}{0.6cm}
\setlength{\leftmargin}{1cm}\setlength{\labelsep}{0.25cm}
\setlength{\parsep}{-3pt}}}{\end{list}~\\[-6ex]}

\newcounter{punkt2}
\newenvironment{enumbib}{\setcounter{punkt2}{0}
\begin{list}{\arabic{punkt2}.}{\usecounter{punkt2}
\setlength{\topsep}{1ex}\setlength{\labelwidth}{0.6cm}
\setlength{\leftmargin}{1cm}\setlength{\labelsep}{0.25cm}
\setlength{\parsep}{1pt}}}{\end{list}}

\newcounter{pic}
\setcounter{pic}{0}

\newcommand{\bdpm}{\begin{displaymath}}
\newcommand{\edpm}{\end{displaymath}}

\newcommand{\beas}{\begin{eqnarray*}}
\newcommand{\eeas}{\end{eqnarray*}}

\newcommand{\ba}{\begin{array}}
\newcommand{\ea}{\end{array}}

\newenvironment{bew}{\vspace*{-0.25cm}\begin{sloppypar}\noindent{\it 
Proof.}}{\hfill\qed\end{sloppypar}\vspace*{0.15cm}}

\newtheorem{satz}{Theorem}[section]          
\newtheorem{lem}[satz]{Lemma}
\newtheorem{kor}[satz]{Corollary}
\newtheorem{prop}[satz]{Proposition}
\newtheorem{bsp}[satz]{Example}
\newtheorem{bem}[satz]{Remark}
\newtheorem{beme}[satz]{Remarks}

\newcommand{\brm}{\begin{rm}}
\newcommand{\erm}{\end{rm}}

\newcommand{\qed}{\hfill $\Box$}
\newcommand{\des}{{\sf des}}
\newcommand{\exc}{{\sf exc}}
\newcommand{\inv}{{\sf inv}}
\newcommand{\rg}{{\rm rank}}

\newcommand{\Ns}
{\unitlength=0.25cm
\begin{picture}(1,1)
\bezier{100}(0,0)(0.5,0.5)(1,1)
\end{picture}}
\newcommand{\Ss}
{\unitlength=0.25cm
\begin{picture}(1,1)
\bezier{100}(0,1)(0.5,0.5)(1,0)
\end{picture}}


\begin{center}
{\large\bf ON THE DIAGRAM OF 132-AVOIDING PERMUTATIONS}\\[1cm]
Astrid Reifegerste\\
Institut f\"ur Mathematik, Universit\"at Hannover\\
Welfengarten 1\\ 
D-30167 Hannover, Germany\\
{\it reifegerste@math.uni-hannover.de}\\[0.5cm]
October 15, 2002 
\end{center}
\vspace*{0.5cm}

\begin{footnotesize}
{\sc Abstract.} The diagram of a $132$-avoiding permutation can easily be characterized: it is 
simply the diagram of a partition. Based on this fact, we present a new bijection between 
$132$-avoiding and $321$-avoiding permutations. We will show that this bijection 
translates the correspondences between these permutations and Dyck paths 
given by Krattenthaler and by Billey-Jockusch-Stanley, respectively, to each 
other. Moreover, the diagram approach yields simple proofs for some enumerative 
results concerning forbidden patterns in $132$-avoiding permutations.
\end{footnotesize}
\vspace*{1cm}


\setcounter{section}{1}\setcounter{satz}{0}

\centerline{\large{\bf 1}\hspace*{0.25cm}
{\sc Introduction}}
\vspace*{0.5cm}

Let ${\cal S}_n$ denote the symmetric group on $\{1,\ldots,n\}$. Given a 
permutation $\pi=\pi_1\cdots\pi_n\in{\cal S}_n$ and a permutation 
$\tau=\tau_1\cdots\tau_k\in{\cal S}_k$, we say that $\pi$ {\it contains the 
pattern} $\tau$ if there is a sequence $1\le i_1<i_2<\ldots<i_k\le n$ such that 
the elements $\pi_{i_1}\pi_{i_2}\cdots\pi_{i_k}$ are in the same relative order 
as $\tau_1\tau_2\cdots\tau_k$. Otherwise, $\pi$ {\it avoids the pattern} 
$\tau$, or alternatively, $\pi$ is {\it $\tau$-avoiding}. We denote by ${\cal 
S}_n(\tau)$ the set of all permutations in ${\cal S}_n$ which avoid $\tau$.\\
It is an often quoted fact that $|{\cal S}_n(\tau)|$ is equal to the $n$th 
Catalan number $C_n=\frac{1}{n+1}{2n\choose n}$ for each pattern $\tau\in{\cal 
S}_3$. Because of obvious symmetry arguments, from an enumerative viewpoint there are only two distinct cases 
to consider, $\tau\in\{123,321\}$ and $\tau\in\{132,213,231,312\}$. Several 
authors established bijections between permutations avoiding a pattern of each 
of these classes. The first one was given by Simion and Schmidt 
\cite{simion-schmidt}; West described in \cite{west} a construction using 
trees; and recently, Krattenthaler \cite{krattenthaler} connected the $123$-avoiding and 
$132$-avoiding permutations via Dyck paths.\\
In Section 2, we present a bijection between ${\cal S}_n(321)$ and ${\cal S}_n(132)$ based on 
another interesting combinatorial object, the diagrams.\\[2ex]
Our construction has the advantage that the excedances of a permutation in ${\cal S}_n(321)$ are 
precisely the descents of its image in ${\cal S}_n(132)$.\\
An {\it excedance of $\pi$} is an integer $i\in\{1,\ldots,n-1\}$ such that 
$\pi_i>i$. Here the element $\pi_i$ is called an {\it excedance letter for $\pi$}. 
Given a permutation $\pi$, we denote the set of excedances of $\pi$ by ${\sf E}(\pi)$ 
and the number $|{\sf E}(\pi)|$ by $\exc(\pi)$. An integer $i\in\{1,\ldots,n-1\}$ for which $\pi_i>\pi_{i+1}$ is called a {\it descent of 
$\pi$}. If $i$ is a descent, we say that $\pi_{i+1}$ is a {\it descent bottom for $\pi$}. 
The set of descents of $\pi$ is denoted by ${\sf D}(\pi)$, its cardinality is denoted by $\des(\pi)$, as usual.\\
Given a permutation $\pi\in{\cal S}_n$, let $i_1<i_2<\ldots<i_e$ be the 
excedances of $\pi$ and $j_1<j_2<\ldots<j_{n-e}$ the remaining positions.
It characterizes $321$-avoiding permutations that both subwords 
$\pi_{i_1}\pi_{i_2}\cdots\pi_{i_e}$ and $\pi_{j_1}\pi_{j_2}\cdots\pi_{j_{n-e}}$ are increasing. 
Thus any $321$-avoiding permutation is uniquely determined by its excedances and 
excedance letters.\\[2ex]
There are several correspondences between restricted permutations and 
lattice paths, in particular, Dyck paths. A {\it Dyck path} is a path in the 
$(x,y)$-plane from the origin to $(2n,0)$ with steps $[1,1]$ (called up-steps) and 
$[1,-1]$ (called down-steps) that never falls below the $x$-axis.\\
For $321$-avoiding permutations such a bijection was given by Billey, Jockusch 
and Stanley \cite{billey etal}; for $132$-avoiding permutations Krattenthaler 
proposed a correspondence to Dyck paths in \cite{krattenthaler}. In Section 3 
we will show that the Dyck path obtained for any $\pi\in{\cal S}_n(321)$ by the 
first mentioned correspondence and the Dyck path associating by Krattenthaler's 
correspondence with the image (with respect to our bijection) $\sigma\in{\cal S}_n(132)$ 
of $\pi$ are the same.\\
Moreover, it will turn out that the diagram of a $132$-avoiding permutation is 
closed related to the corresponding Dyck path. This yields a 
simple explanation for the connections between the number of inversions of the 
permutation and several parameters of the Dyck path recently appeared in 
\cite{braenden etal}. Further, we use this relation to enumerate 
restricted partitions of prescribed rank.\\[2ex]
In Section 4 the diagram approach will be used to obtain some enumerative 
results concerning the restriction of $132$-avoiding permutations by additional 
patterns. These results are already known (see \cite{mansour-vainshtein1}) but we will give 
a bijective proof of them.\\
The paper ends with a note on how to obtain the number of occurrences of the pattern 
$132$ in an arbitrary permutation via the diagram.
\vspace*{0.75cm}


\setcounter{section}{2}\setcounter{satz}{0}

\centerline{\large{\bf 2}\hspace*{0.25cm}
{\sc A bijection between $132$-avoiding and $321$-avoiding permutations}}
\vspace*{0.5cm}

Let ${\cal Y}_n:=\{(\lambda_1,\ldots,\lambda_{n-1}):0\le\lambda_{n-1}\le
\lambda_{n-2}\le\ldots\le\lambda_1\le n-1,\;\lambda_i\le n-i\mbox{ for all }i\}$ be the set of 
partitions whose Young diagram fits in the shape $(n-1,n-2,\ldots,1)$. (We will identify a partition with its Young diagram and 
vice versa.) In \cite[Sect. 3.2]{reifegerste} (or \cite{reifegerste1}), we 
describe a bijection ${\cal S}_n(321)\to{\cal Y}_n$
which takes the permutation $\pi$ with excedances $i_1,\ldots,i_e$ to the diagram  
with corners $(i_k,n+1-\pi_{i_k})$ where $k=1,\ldots,e$. For $132$-avoiding permutations a simple correspondence to 
partitions with restricted diagram can be given, as well.\\[2ex] 
The key object in our derivation is the {\it diagram of a permutation} (for an introduction see 
\cite[chap. 1]{macdonald}). Given a permutation $\pi\in{\cal S}_n$, we obtain the diagram $D(\pi)$ of $\pi$ as 
follows. Let $\pi$ be represented by an $n\times n$-array with a dot in 
each of the squares $(i,\pi_i)$. (The other cells are white.) Shadow all squares due south or due east of some 
dot and the dotted cell itself. The diagram $D(\pi)$ is defined as the region 
left unshaded after this procedure.
 
\begin{bsp} \label{Beispiel 1}
\brm
The diagram of $\pi=4\:2\:8\:3\:6\:9\:7\:5\:1\:10\in{\cal S}_{10}$ contains the 
white squares of
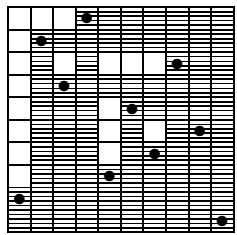
\begin{center}
\unitlength=0.3cm
\begin{picture}(10,10)
\linethickness{0.3pt}
\multiput(0,0)(0,1){11}{\line(1,0){10}}
\multiput(0,0)(1,0){11}{\line(0,1){10}}
\put(3.5,9.5){\circle*{0.5}}
\put(1.5,8.5){\circle*{0.5}}
\put(7.5,7.5){\circle*{0.5}}
\put(2.5,6.5){\circle*{0.5}}
\put(5.5,5.5){\circle*{0.5}}
\put(8.5,4.5){\circle*{0.5}}
\put(6.5,3.5){\circle*{0.5}}
\put(4.5,2.5){\circle*{0.5}}
\put(0.5,1.5){\circle*{0.5}}
\put(9.5,0.5){\circle*{0.5}}
\linethickness{0.2pt}
\multiput(3,9)(0,0.2){5}{\line(1,0){7}}
\multiput(1,8)(0,0.2){5}{\line(1,0){9}}
\multiput(7,7)(0,0.2){5}{\line(1,0){3}}
\multiput(2,6)(0,0.2){5}{\line(1,0){8}}
\multiput(5,5)(0,0.2){5}{\line(1,0){5}}
\multiput(8,4)(0,0.2){5}{\line(1,0){2}}
\multiput(6,3)(0,0.2){5}{\line(1,0){4}}
\multiput(4,2)(0,0.2){5}{\line(1,0){6}}
\multiput(0,1)(0,0.2){5}{\line(1,0){10}}
\multiput(0,0)(0,0.2){5}{\line(1,0){10}}
\multiput(1,2)(0,0.2){5}{\line(1,0){3}}
\multiput(1,3)(0,0.2){5}{\line(1,0){3}}
\multiput(5,3)(0,0.2){5}{\line(1,0){1}}
\multiput(1,4)(0,0.2){5}{\line(1,0){3}}
\multiput(1,5)(0,0.2){5}{\line(1,0){3}}
\multiput(1,6)(0,0.2){5}{\line(1,0){1}}
\multiput(1,7)(0,0.2){5}{\line(1,0){1}}
\multiput(3,7)(0,0.2){5}{\line(1,0){1}}
\multiput(5,4)(0,0.2){5}{\line(1,0){1}}
\multiput(7,4)(0,0.2){5}{\line(1,0){1}}
\end{picture}
\end{center}
\erm
\end{bsp}  

By the construction, each of the connected components of $D(\pi)$ is a Young diagram. 
Their corners are defined to be the elements of the {\it essential set} ${\cal E}(\pi)$ of the 
permutation $\pi$. In \cite{fulton}, Fulton introduced this set which together 
with a rank function was used as a tool for algebraic treatment of Schubert 
polynomials. In \cite{eriksson-linusson}, Eriksson and Linusson characterized the essential sets that can 
arise from arbitrary permutations, as well as those coming from certain classes of 
permutations.\\[2ex]
It is very easy to characterize the diagrams of $132$-avoiding permutations. 

\begin{satz} \label{132-avoiding diagram}
Let $\pi\in{\cal S}_n$ be a permutation not equal to the identity. Then $\pi$ is 
$132$-avoiding if and only if its diagram consists of only one component and 
$(1,1)\in D(\pi)$.
\end{satz} 

\begin{bew}
If there are indices $i<j<k$ such that $\pi_i<\pi_k<\pi_j$, then the square 
$(j,\pi_k)$ belongs to $D(\pi)$, but it is not connected with $(1,1)$:
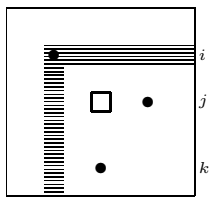
\begin{center}
\unitlength=0.25cm
\begin{picture}(10,10)
\linethickness{0.2pt}
\multiput(0,0)(0,10){2}{\line(1,0){10}}
\multiput(0,0)(10,0){2}{\line(0,1){10}}
\put(2.5,7.5){\circle*{0.5}}
\put(5,1.5){\circle*{0.5}}
\put(7.5,5){\circle*{0.5}}
\multiput(2,7)(0,0.2){6}{\line(1,0){8}}
\multiput(2,0)(0,0.2){35}{\line(1,0){1}}
\put(10.2,7.5){\makebox(0,0)[lc]{\tiny$i$}}
\put(10.2,5){\makebox(0,0)[lc]{\tiny$j$}}
\put(10.2,1.5){\makebox(0,0)[lc]{\tiny$k$}}
\linethickness{0.7pt}
\put(4.5,5.5){\line(1,0){1}}\put(4.5,4.5){\line(1,0){1}}
\put(4.5,4.5){\line(0,1){1}}\put(5.5,4.5){\line(0,1){1}}
\end{picture}
\end{center}
Clearly, the existence of such a square is also sufficient for $\pi$ containing the 
pattern $132$.\\ 
Note that the square $(1,1)$ must be an element of $D(\pi)$ for any $132$-avoiding permutation $\pi\not={\rm id}$, 
otherwise we would have $\pi_1=1$ and hence $\pi_i=i$ for all $i=1,\ldots,n$.
\end{bew}

Thus the diagram $D(\pi)$ of a permutation $\pi\in{\cal S}_n(132)$ 
is the graphical representation of a partition. By construction, $D(\pi)$ is 
just the diagram of an element of ${\cal Y}_n$: the square $(i,j(i))$ belongs 
to $D(\pi)$ if and only if no index $k\le i$ satisfies $\pi_k\le j$. Thus we have $j(i)\le n-i$.\\[2ex]
This yields a simple bijection between $321$-avoiding and $132$-avoiding 
permutations on $\{1,\ldots,n\}$ which is in the following denoted by $\Phi$.\\ 
Let $i_1,\ldots,i_e$ be the excedances of $\pi\in{\cal S}_n(321)$. Then the diagram 
of the corresponding permutation $\Phi(\pi)$ is equal to the Young diagram with corners 
$(i_k,n+1-\pi_{i_k}),\;k=1,\ldots,e$. 
Recovering $\Phi(\pi)$ from $D(\Phi(\pi))$ is trivial: row by row, put a dot in the leftmost shaded square such that there is exactly one dot in each column.   

\begin{bsp}
\brm
For the permutation $\pi=1\:4\:7\:2\:3\:8\:5\:6\:10\:9\in{\cal S}_{10}(321)$ we 
have ${\sf E}(\pi)=\{2,3,6,9\}$. Hence it corresponds to the permutation with the diagram
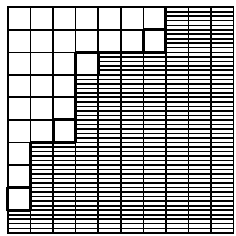
\begin{center}
\unitlength=0.3cm
\begin{picture}(10,10)
\linethickness{0.2pt}
\multiput(0,0)(0,1){11}{\line(1,0){10}}
\multiput(0,0)(1,0){11}{\line(0,1){10}}
\multiput(7,9)(0,0.2){5}{\line(1,0){3}}
\multiput(7,8)(0,0.2){5}{\line(1,0){3}}
\multiput(4,7)(0,0.2){5}{\line(1,0){6}}
\multiput(3,6)(0,0.2){5}{\line(1,0){7}}
\multiput(3,5)(0,0.2){5}{\line(1,0){7}}
\multiput(3,4)(0,0.2){5}{\line(1,0){7}}
\multiput(1,3)(0,0.2){5}{\line(1,0){9}}
\multiput(1,2)(0,0.2){5}{\line(1,0){9}}
\multiput(1,1)(0,0.2){5}{\line(1,0){9}}
\multiput(0,0)(0,0.2){5}{\line(1,0){10}}
\linethickness{0.7pt}
\put(4,8){\line(1,0){3}}\put(6,9){\line(1,0){1}}
\put(6,8){\line(0,1){1}}\put(7,8){\line(0,1){2}}
\put(3,8){\line(1,0){1}}\put(3,7){\line(1,0){1}}
\put(3,4){\line(0,1){4}}\put(4,7){\line(0,1){1}}
\put(2,5){\line(1,0){1}}\put(1,4){\line(1,0){2}}
\put(2,4){\line(0,1){1}}\put(1,1){\line(0,1){3}}
\put(0,1){\line(1,0){1}}\put(0,2){\line(1,0){1}}
\put(0,1){\line(0,1){1}}
\end{picture}
\end{center}
and we obtain $\Phi(\pi)=8\:9\:5\:4\:6\:7\:2\:3\:10\:1\in{\cal S}_{10}(132)$. 
\erm
\end{bsp}  

As observed by Fulton in \cite{fulton}, every row of a permutation diagram 
containing a white corner (that is an element of the essential set) corresponds to a descent. Thus we 
have $\des(\Phi(\pi))=\exc(\pi)$ for all $\pi\in{\cal S}_n(321)$. But there is 
more to it than that: the excedance set of $\pi$ and the descent set of 
$\Phi(\pi)$ have not only the same number of elements; 
the sets are even identical. 

\begin{prop}
We have ${\sf E}(\pi)={\sf D}(\Phi(\pi))$ for all $\pi\in{\cal S}_n(321)$.
\end{prop}

\begin{bew}
Any excedance $i$ of $\pi$ corresponds to a corner $(i,n+1-\pi_i)$ of $D(\Phi(\pi))$. 
Obviously, by constructing $\Phi(\pi)$ from its diagram we obtain a 
descent of $\Phi(\pi)$ at the position $i$. 
\end{bew}

Every $321$-avoiding permutation is completely determined by its excedances and 
excedance letters. Our bijection shows that it is sufficient for fixing a $132$-avoiding 
permutation to know the descents, the descent bottoms, and the first letter.\\
Let $i_1<\ldots<i_e$ be the excedances of $\pi\in{\cal S}_n(321)$, and let $\sigma:=\Phi(\pi)$. 
Then we have 
\beas
&&\sigma_1=n+2-\pi_{i_1},\\
&&\sigma_{i_k+1}=n+2-\pi_{i_{k+1}}\quad\mbox{for }k=1,\ldots,e-1,\\
&&\sigma_{i_e+1}=1.
\eeas           
It is clear from the construction that these elements are precisely the left-to-right minima 
of $\sigma$. (A {\it left-to-right minimum} of a permutation $\sigma$ is an element 
$\sigma_i$ which is smaller than all elements to its left, i.e., 
$\sigma_i<\sigma_j$ for every $j<i$.) Based on this, we can determine 
the permutation $\sigma$ since $\sigma$ is $132$-avoiding.

\begin{bsp}
\brm
Let again $\pi=1\:\underline{4}\:\underline{7}\:2\:3\:\underline{8}\:5\:6\:\underline{10}\:9
\in{\cal S}_{10}(321)$. (The underlined positions are just the excedances of $\pi$.) As 
described above, we obtain the left-to-right minima of $\Phi(\pi)$ and their 
positions: 
\bdpm
8\:\ast\:5\:4\:\ast\:\ast\:2\:\ast\:\ast\:1,
\edpm 
and hence, by putting the remaining elements $a=3,6,7,9,10$ on the first possible position 
following $a-1$, the permutation $\Phi(\pi)=8\:9\:5\:4\:6\:7\:2\:3\:10\:1$. 
\erm
\end{bsp}

\begin{beme} \label{remarks to section 2}
\brm
\begin{enum}
\item[]
\item In Chapter 1 of \cite{macdonald}, Macdonald defined the {\it 
dominant} permutations. This special case of {\it vexillary} (or $2143$-avoiding) permutations 
is characterized by the following equivalent conditions:
\begin{enum2}
\item $D(\pi)$ is the diagram of a partition.
\item The code $c$ of $\pi$ is a partition, that is, $c_i\ge c_{i+1}$ for all 
$i=1,\ldots,n-1$.
\end{enum2}
\vspace*{-0.2cm}

(For a permutation $\pi$ the $i$th component of its {\it code} counts the 
number of indices $j>i$ satisfying $\pi_j<\pi_i$.) As mentioned above, we may 
assume in (i) that $D(\pi)$ is an element of ${\cal Y}_n$. Thus ${\cal 
S}_n(132)$ is precisely the set of all dominant permutations. In particular, there 
are $C_n$ such permutations.
\item For {\it any} permutation $\pi\in{\cal S}_n$ the number of squares in the $i$th row of 
$D(\pi)$ is equal to the $i$th component $c_i$ of its code (see 
\cite[p. 9]{macdonald}). Hence it follows that 
$|D(\pi)|=c_1+\ldots+c_n=\inv(\pi)$, where $\inv(\pi)$ denotes the number of 
inversions of $\pi$.
\item In \cite[Sect. 3.2]{reifegerste} (or in \cite{reifegerste1}), we prove that the number of excedances 
has the Narayana distribution over ${\cal S}_n(321)$. Consequently, the number 
of diagrams fitting in the shape $(n-1,n-2,\ldots,1)$ and having $k$ corners is given by the Narayana 
number $N(n,k+1)=\frac{1}{n}{n\choose k}{n\choose k+1}$. Thus the bijection immediately 
shows that the statistic $\des$ is Narayana distributed over ${\cal S}_n(132)$.
\item Also in \cite[Sect. 2.3]{reifegerste}, an involution on ${\cal S}_n(321)$ 
was established which proves the symmetry of the joint distribution of the pair 
$(\exc,\inv)$ over ${\cal S}_n(321)$. For $\pi\in{\cal S}_n(321)$ let $\sigma=\sigma_1\cdots\sigma_n$ be the image 
with respect to this map. Then the reverse $\sigma_n\cdots\sigma_1\in{\cal S}_n(123)$ 
corresponds to $\Phi(\pi)\in{\cal S}_n(132)$ by the bijection due to Simion 
and Schmidt (see \cite[Prop. 19]{simion-schmidt}).
\item In \cite[p. 7]{fulmek}, Fulmek has given a "`graphical"' construction of 
a map from ${\cal S}_n$ to the set of Dyck paths of length $2n$ whose restrictions 
on ${\cal S}_n(321)$ and ${\cal S}_n(312)$, respectively, are bijections. 
Based on the graph of a permutation $\pi\in{\cal S}_n$ he considered the union 
of all points which are in the southeast of some left-to-right maximum. 
The corresponding path is defined to be the boundary of this region. It is easy 
to see that the region above the path is just the connected component 
of all squares of rank zero belonging to the diagram of 
the inverse of $\pi_n\pi_{n-1}\cdots\pi_1$. (See also \cite{eriksson-linusson} for a description how to obtain an arbitrary permutation from its ranked essential set.)   
\end{enum}                
\erm
\end{beme}
\vspace*{0.5cm}


\setcounter{section}{3}\setcounter{satz}{0}

\centerline{\large{\bf 3}\hspace*{0.25cm}
{\sc Correspondences to Dyck Paths}}
\vspace*{0.5cm}

For $321$-avoiding, as well as for $132$-avoiding permutations one-to-one correspondences 
to lattice paths were given by several authors. In \cite[p. 361]{billey etal}, Billey, 
Jockusch and Stanley established a bijection $\Psi_{BJS}$ between $321$-avoiding 
permutations on $\{1,\ldots,n\}$ and Dyck paths of length $2n$. Recently in 
\cite[Sect. 2]{krattenthaler}, Krattenthaler 
exhibited a Dyck path correspondence $\Psi_K$ for $132$-avoiding permutations. Our bijection $\Phi$ 
translates these constructions into each other.

\begin{satz}
Let $\pi\in{\cal S}_n(321)$. Then we have $\Psi_{BJS}(\pi)=\Psi_{K}(\Phi(\pi))$.
\end{satz}

\begin{bew}
Let $\pi\in{\cal S}_n(321)$ with the excedances $i_1<\ldots<i_e$, and let 
$\sigma:=\Phi(\pi)$. The bijection $\Psi_{BJS}$ constructs the Dyck path corresponding to $\pi$ as follows:
\begin{enum1} 
\item Let $a_k:=\pi_{i_k}-1$ for $k=1,\ldots,e$ and $a_0:=0,\;a_{e+1}:=n$. 
Furthermore, let $b_k:=i_k$ for $k=1,\ldots,e$ and $b_0:=0,\;b_{e+1}:=n$. 
\item Generate the Dyck path (starting at the origin) by adjoining $a_k-a_{k-1}$ up-steps and $b_k-b_{k-1}$ down-steps 
for $k=1,\ldots,e+1$. 
\end{enum1}
\vspace*{-0.35cm}

As shown in the preceding section, the elements
\beas
c_1&:=&\sigma_1\;=\;n+2-\pi_{i_1},\\
c_{k+1}&:=&\sigma_{i_k+1}\;=\;n+2-\pi_{i_{k+1}}\quad\mbox{for }k=1,\ldots,e-1,\\
c_{e+1}&:=&\sigma_{i_e+1}\;=\;1
\eeas           
are the left-to-right minima of $\sigma$. With the convention $c_0:=n+1$ we 
have $c_{k-1}-c_k=a_k-a_{k-1}$ for all $k=1,\ldots,e+1$. For the number $d_k$ of the 
positions between the $k$th and (including) the $(k+1)$st left-to-right minimum 
we obtain $d_k=b_k-b_{k-1}$ for $k=1,\ldots,e+1$. (Let $n+1$ be the position of 
the imaginary $(e+2)$nd minimum, so $d_{e+1}=n-b_e$.)\\
Hence the translation of $\Psi_{BJS}$ by $\Phi$ constructs the Dyck path corresponding to $\sigma\in{\cal S}_n(132)$ as follows:
\begin{enum1} 
\item Let $c_1>\ldots>c_{e+1}$ be the left-to-right minima of $\sigma$. 
Furthermore, let $d_k$ be one plus the number of the letters in $\sigma$ between $c_k$ and $c_{k+1}$ for $k=1,\ldots,e+1$. Initialize $c_0:=n+1$.
\item Generate the Dyck path (starting at the origin) by adjoining $c_{k-1}-c_k$ up-steps and $d_k$ down-steps 
for $k=1,\ldots,e+1$. 
\end{enum1}
\vspace*{-0.35cm}

But this is precisely the description of $\Psi_{K}$ proposed in \cite{krattenthaler}.
\end{bew}

\begin{bsp}
\brm
Let $\pi=1\:4\:7\:2\:3\:8\:5\:6\:10\:9\in{\cal S}_{10}(321)$, and let 
$\sigma=\Phi(\pi)=8\:9\:5\:4\:6\:7\:2\:3\:10\:1$. 
Billey-Jockusch-Stanley's bijection takes $\pi$ to the Dyck path
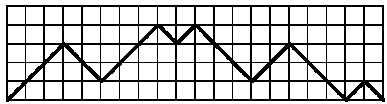
\begin{center}
\unitlength=0.25cm
\begin{picture}(20,5)
\linethickness{0.2pt}
\multiput(0,0)(0,1){6}{\line(1,0){20}}
\multiput(0,0)(1,0){21}{\line(0,1){5}}
\linethickness{0.7pt}
\bezier{300}(0,0)(1.5,1.5)(3,3)
\bezier{200}(3,3)(4,2)(5,1)
\bezier{300}(5,1)(6.5,2.5)(8,4)
\bezier{100}(8,4)(8.5,3.5)(9,3)
\bezier{100}(9,3)(9.5,3.5)(10,4)
\bezier{300}(10,4)(11.5,2.5)(13,1)
\bezier{200}(13,1)(14,2)(15,3)
\bezier{300}(15,3)(16.5,1.5)(18,0)
\bezier{100}(18,0)(18.5,0.5)(19,1)
\bezier{100}(19,1)(19.5,0.5)(20,0)
\end{picture}
\end{center}
which is exactly the path corresponding to $\sigma$ by Krattenthaler's bijection.
\erm
\end{bsp}  

The following results use the (now obvious) fact that the Dyck path $\Psi_K(\pi)$ and the 
diagram of a $132$-avoiding permutation $\pi$ are closely related to each 
other. Given a permutation $\pi\in{\cal S}_n(132)$, its diagram $D(\pi)$ is 
just the region bordered by the lines between the lattice points $(0,0)$ and 
$(n,n)$ and between $(n,n)$ and $(2n,0)$, respectively, and the path $\Psi_K(\pi)$. (The northwest-to-southeast diagonals correspond 
to the diagram columns.)

\begin{bsp}
\brm
For $\pi=8\:9\:5\:4\:6\:7\:2\:3\:10\:1\in{\cal S}_{10}(132)$ we obtain:
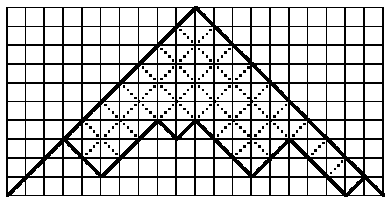
\begin{center}
\unitlength=0.25cm
\begin{picture}(20,10)
\linethickness{0.2pt}
\multiput(0,0)(0,1){11}{\line(1,0){20}}
\multiput(0,0)(1,0){21}{\line(0,1){10}}
\linethickness{0.7pt}
\bezier{1000}(0,0)(5,5)(10,10)
\bezier{200}(3,3)(4,2)(5,1)
\bezier{300}(5,1)(6.5,2.5)(8,4)
\bezier{100}(8,4)(8.5,3.5)(9,3)
\bezier{100}(9,3)(9.5,3.5)(10,4)
\bezier{300}(10,4)(11.5,2.5)(13,1)
\bezier{200}(13,1)(14,2)(15,3)
\bezier{300}(15,3)(16.5,1.5)(18,0)
\bezier{100}(18,0)(18.5,0.5)(19,1)
\bezier{1000}(10,10)(15,5)(20,0)
\linethickness{0.5pt}
\bezier{35}(4,2)(7.5,5.5)(11,9)
\bezier{20}(8,4)(10,6)(12,8)
\multiput(10,4)(1,-1){3}{\bezier{15}(0,0)(1.5,1.5)(3,3)}
\multiput(15,3)(1,-1){3}{\bezier{5}(0,0)(0.5,0.5)(1,1)}
\bezier{40}(9,9)(13,5)(17,1)
\bezier{30}(8,8)(11,5)(14,2)
\bezier{15}(7,7)(8.5,5.5)(10,4)
\multiput(4,2)(1,1){3}{\bezier{10}(0,2)(1,1)(2,0)}
\end{picture}
\end{center}
\erm
\end{bsp}  

In \cite{braenden etal}, the authors studied the statistic $e_k$ that counts the 
number of increasing subsequences of length $k+1$ in a permutation. Expressing 
$e_k(\pi)$ in terms of the Dyck path $\Psi_K(\pi)$, some 
applications were given in \cite{braenden etal} which relate various combinatorial structures to 
$132$-avoiding permutations. The translation of the statistic into Dyck 
path characteristics becomes now immediately clear.

\begin{kor}
For any permutation $\pi\in{\cal S}_n(132)$ consider the associated Dyck path 
$\Psi_K(\pi)$ and denote the height of the starting point of its $i${\rm th} step by $w_i(\pi)$. Then we have
\begin{enum}
\item $w_1(\pi)+\ldots+w_{2n}(\pi)=n^2-2\:\inv(\pi)$.
\item $w_{i_1}(\pi)+\ldots+w_{i_s}(\pi)={n+1\choose 2}-\inv(\pi)$, where 
$i_1,\ldots,i_s$ are the indices of the down-steps.
\end{enum}
\end{kor}
\vspace*{-0.45cm}

\begin{bew}
As remarked in \ref{remarks to section 2}b), we have $|D(\pi)|=\inv(\pi)$ for all $\pi\in{\cal S}_n$.\\
{\bf a)} The sum of heights of the Dyck path $(\Ns)^n(\Ss)^n$ (where exponentiation 
denotes repetition) is equal to $n^2$. Any 
square of $D(\pi)$ reduces this value by 2:
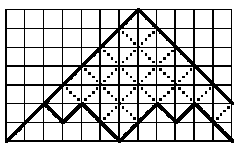
\begin{center}
\unitlength=0.25cm
\begin{picture}(12,7)
\linethickness{0.2pt}
\multiput(0,0)(0,1){8}{\line(1,0){12}}
\multiput(0,0)(1,0){13}{\line(0,1){7}}
\linethickness{0.7pt}
\bezier{700}(0,0)(3.5,3.5)(7,7)
\bezier{100}(2,2)(2.5,1.5)(3,1)
\bezier{100}(3,1)(3.5,1.5)(4,2)
\bezier{200}(4,2)(5,1)(6,0)
\bezier{200}(6,0)(7,1)(8,2)
\bezier{100}(8,2)(8.5,1.5)(9,1)
\bezier{100}(9,1)(9.5,1.5)(10,2)
\bezier{200}(10,2)(11,1)(12,0)
\bezier{500}(12,2)(9.5,4.5)(7,7)
\linethickness{0.5pt}
\bezier{20}(4,2)(6,4)(8,6)
\bezier{20}(5,1)(7,3)(9,5)
\bezier{10}(8,2)(9,3)(10,4)
\bezier{5}(10,2)(10.5,2.5)(11,3)
\bezier{5}(11,1)(11.5,1.5)(12,2)
\bezier{5}(3,3)(3.5,2.5)(4,2)
\bezier{15}(4,4)(5.5,2.5)(7,1)
\bezier{15}(5,5)(6.5,3.5)(8,2)
\bezier{20}(6,6)(8,4)(10,2)
\end{picture}
\end{center}
{\bf b)} The sum of all down-step heights counts the number of squares of the 
slanting lattices below the path:  
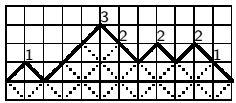
\begin{center}
\unitlength=0.25cm
\begin{picture}(12,5)
\linethickness{0.2pt}
\multiput(0,0)(0,1){6}{\line(1,0){12}}
\multiput(0,0)(1,0){13}{\line(0,1){5}}
\linethickness{0.7pt}
\put(0,1){\line(1,0){12}}
\bezier{100}(0,1)(0.5,1.5)(1,2)
\bezier{100}(1,2)(1.5,1.5)(2,1)
\bezier{300}(2,1)(3.5,2.5)(5,4)
\bezier{200}(5,4)(6,3)(7,2)
\bezier{100}(7,2)(7.5,2.5)(8,3)
\bezier{100}(8,3)(8.5,2.5)(9,2)
\bezier{100}(9,2)(9.5,2.5)(10,3)
\bezier{200}(10,3)(11,2)(12,1)
\linethickness{0.5pt}
\bezier{5}(0,1)(0.5,0.5)(1,0)
\bezier{5}(1,0)(1.5,0.5)(2,1)
\bezier{5}(2,1)(2.5,0.5)(3,0)
\bezier{15}(3,0)(4.5,1.5)(6,3)
\multiput(5,0)(2,0){3}{\bezier{10}(0,0)(1,1)(2,2)}
\bezier{5}(11,0)(11.5,0.5)(12,1)
\bezier{10}(3,2)(4,1)(5,0)
\bezier{15}(4,3)(5.5,1.5)(7,0)
\multiput(7,0)(2,0){2}{\bezier{10}(0,2)(1,1)(2,0)}
\put(1.25,2.4){\makebox(0,0)[cc]{\tiny\sf 1}}
\put(5.25,4.4){\makebox(0,0)[cc]{\tiny\sf 3}}
\put(6.25,3.4){\makebox(0,0)[cc]{\tiny\sf 2}}
\put(8.25,3.4){\makebox(0,0)[cc]{\tiny\sf 2}}
\put(10.25,3.4){\makebox(0,0)[cc]{\tiny\sf 2}}
\put(11.25,2.4){\makebox(0,0)[cc]{\tiny\sf 1}}
\end{picture}
\end{center}
For the path $(\Ns)^n(\Ss)^n$ this number is ${n+1\choose 2}$.
\end{bew}

\begin{beme}
\brm
\begin{enum}
\item[]
\item The sum of the heights of all steps is easily seen to equal the area of the 
Dyck path. It was shown in \cite[Sect. 3.1]{braenden etal} that the sum relates 
to the statistic $2e_1+e_0$ for $132$-avoiding permutations.
\item For a correspondence between fountains of coins and $132$-avoiding 
permutations, \cite{braenden etal} used that $e_1(\pi)+e_0(\pi)$ is equal to the sum of 
the heights of the down-steps in $\Psi_K(\pi)$.
\end{enum}
\erm
\end{beme}

On the other hand, from \cite{braenden etal} we can derive enumerative results for the partitions in ${\cal 
Y}_n$ and (using the bijection $\Phi$) for $321$-avoiding permutations.\\[2ex]
By \cite[Prop. 7]{braenden etal}, the distribution of right-to-left maxima in ${\cal S}_n(132)$ 
is given by the ballot numbers. (A {\it right-to-left maxima} of a permutation 
$\pi$ is an element $\pi_i$ which is larger than all $\pi_j$ with $j>i$.) The 
number of permutations in ${\cal S}_n(132)$ with $k$ such maxima equals the 
ballot number
\bdpm
b(n-1,n-k)=\frac{k}{2n-k}\:{2n-k\choose n}.
\edpm
Among other things, the {\it ballot number} $b(n,k)$ counts the number of the lattice paths from $(0,0)$ to $(n+k,n-k)$ 
with up-steps and down-steps only, never falling below the $x$-axis. It is well known that $b(n,k)={n+k\choose n}-{n+k\choose 
n+1}=\frac{n+1-k}{n+1}{n+k\choose n}$. (For instance, see \cite[p. 73]{feller}.)\\ 
For their bijective proof, \cite{braenden etal} used Krattenthaler's map: $\Psi_K$ translates any right-to-left maximum of $\pi\in{\cal S}_n(132)$
into a return of the associated Dyck path. (A {\it return} of a Dyck path is a 
down-step landing on the $x$-axis.). The number of returns of Dyck paths is known to have a distribution given by 
$b(n-1,n-k)$. (This result was given in \cite{deutsch}.)

\begin{kor}
The number of permutations $\pi\in{\cal S}_n(321)$ with $k-1$ elements $\pi_i=i+1$ 
is equal to $b(n-1,n-k)$.
\end{kor}

\begin{bew}
By the construction of $\Psi_{BJS}(\pi)$, a return, except the last (down-) 
step, appears if and only if $i$ is an excedance of $\pi$ with $\pi_i=i+1$. The 
very last step of $\Psi_{BJS}(\pi)$ is a return by definition. 
\end{bew}

\begin{kor}
The number of diagrams fitting in $(n-1,n-2,\ldots,1)$ with $k-1$ corners in the diagonal $i+j=n$ equals $b(n-1,n-k)$.
\end{kor}

\begin{bem}
\brm
In comparison with this, the Narayana number $N(n,k)=\frac{1}{n}{n\choose k-1}{n\choose 
k}$ counts the number of restricted diagrams with $k-1$ corners (see also Remark 
\ref{remarks to section 2}c)).\\
The condition $i+j=n$ on a corner $(i,j)$ of the diagram $D(\pi)$ also occurs in the 
following section in context with the pattern $213$ in a $132$-avoiding permutation $\pi$.
\erm
\end{bem}

We can also use the relation between diagrams and Dyck paths to obtain more 
information about the restricted partitions. For instance, the number of 
partitions in ${\cal Y}_n$ of a prescribed rank can easily be derived from a known 
result for paths. (The {\it (Durfee) rank} of a partition $\lambda$, 
denoted by $\rg(\lambda)$, is the largest integer $i$ for which $\lambda_i\ge 
i$, or equivalently, the length of the main diagonal of the diagram of 
$\lambda$.)

\begin{lem}
Let $\pi\in{\cal S}_n(132)$. With the above notation we have 
$w_{n+1}(\pi)=n-2\cdot\rg(D(\pi))$.
\end{lem}

\begin{bew}
This is an immediate consequence of the relation between $D(\pi)$ and 
$\Psi_K(\pi)$.
\end{bew}

\begin{satz}
Let $r(n,k)$ denote the number of partitions in ${\cal Y}_n$ of rank $k$. Then 
we have
\bdpm
r(n,k)=\left(\frac{n+1-2k}{n+1-k}\:{n\choose k}\right)^2\quad\mbox{for all }
0\le k\le\mbox{$\lfloor\frac{n}{2}\rfloor$}.
\edpm
\end{satz}
\vspace*{0.3cm}

\begin{bew}
Given a partition $\lambda\in{\cal Y}_n$ of rank $k$, let $\pi\in{\cal S}_n(132)$ 
be the (uniquely determined) permutation whose diagram equals the diagram of $\lambda$. By the 
lemma, both the first "half" and the second "half" of $\Psi_K(\pi)$ are paths from $(0,0)$ to $(n,n-2k)$ 
only consisting of up-steps and down-steps and never falling below the $x$-axis (for 
the latter one, consider the reflection on the line $x=n$). Thus, $r(n,k)$ is 
the square of the number of precisely those paths. Feller's Ballot Theorem (see 
\cite[p. 73]{feller}) enumerates the paths from the origin to $(n,t)$ with 
$n,t\ge1$ which neither touch (that is, there is no return) or cross the $x$-axis. 
Their number is equal to $\frac{t}{n}{n\choose\frac{n-t}{2}}$. Consequently, there are 
\bdpm
\frac{t+1}{n+1}\:{n+1\choose\frac{n-t}{2}}
\edpm
paths from $(0,0)$ to $(n,t)$ where $n,t\ge0$ with up-steps and down-steps never going below the 
$x$-axis (returns are allowed). To see this, insert an up-step before the first 
step and define its starting point to be the origin. The resulting path satisfies the conditions of Feller's theorem. 
Conversely, each path from $(0,0)$ to $(n+1,t+1)$ without returns can easily be transformed into 
an admissible path by deleting the first (up-) step and redefining the origin. 
Taking $t=n-2k$ we obtain
\bdpm
r(n,k)=\bigg(\frac{n-2k+1}{n+1}\:{n+1\choose k}\bigg)^2=\bigg(\frac{n+1-2k}{n+1-k}\:{n\choose k}\bigg)^2
\edpm                
for all $k=0,\ldots,\lfloor\frac{n}{2}\rfloor$.
\end{bew}

\begin{bem}
\brm
By the proof, the positive root $q(n,k)$ of $r(n,k)$ counts the number of paths which 
begin at the origin, end at $(n,n-2k)$, consist of 
up-steps and down-steps, and stay above the $x$-axis. Clearly, $q(n,0)=1$ and 
$q(n,\lfloor\frac{n}{2}\rfloor)=C_{\lceil\frac{n}{2}\rceil}$ for all $n$. 
Obviously, $q(n,k)$ satisfies the recurrence $q(n,k)=q(n-1,k-1)+q(n-1,k)$. (The 
right-hand side is just the number of paths ending at $(n-1,n-2k+1)$ and $(n-1,n-2k-1)$, respectively.) 
Thus the numbers $q(n,k)$ are exactly the entries in the counter diagonals of the Catalan's triangle.
\erm
\end{bem}
\vspace*{0.5cm}


\setcounter{section}{4}\setcounter{satz}{0}

\centerline{\large{\bf 4}\hspace*{0.25cm}
{\sc Forbidden patterns in 132-avoiding permutations}}
\vspace*{0.5cm}

Now we will use the correspondence between ${\cal S}_n(132)$ and ${\cal Y}_n$
for the enumeration of multiple restrictions on permutations. The results 
concerning the Wilf-equivalence of several pairs $\{132,\tau\}$ where $\tau\in{\cal S}_k$ 
are already known, see \cite{mansour-vainshtein2}, \cite{mansour-vainshtein3}, 
\cite{krattenthaler}, \cite{chow-west}. (We say that $\{132,\tau_1\}$ and 
$\{132,\tau_2\}$ are {\it Wilf-equivalent} if $|{\cal S}_n(132,\tau_1)|=|{\cal 
S}_n(132,\tau_2)|$ for all $n$.) While the proofs given in these papers 
are analytical we present bijective ones. 

\begin{satz} \label{Diagrammcharakterisierung}
Let $\pi\in{\cal S}_n(132)$ be a permutation not equal to the identity, and 
$k\ge3$. Then 
\begin{enum}
\item $\pi$ avoids $k(k-1)\cdots1$ if and only if $D(\pi)$ has at 
most $k-2$ corners. In particular, we have $\des(\pi)\le k-2$.
\item $\pi$ avoids $12\cdots k$ if and only if $D(\pi)$ contains the diagram 
$(n+1-k,n-k,\ldots,1)$.
\item $\pi$ avoids $213\cdots k$ if and only if every corner $(i,j)$ of $D(\pi)$ 
satisfies $i+j\ge n+3-k$. 
\end{enum}
\end{satz}
\vspace*{-0.45cm}

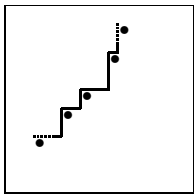
\begin{bew}\\
{\bf a)}\hspace*{0.5cm} 
\begin{minipage}[t]{15cm}
\raisebox{-2.3cm}
{\unitlength=0.25cm
\begin{picture}(10,10)
\linethickness{0.2pt}
\multiput(0,0)(0,10){2}{\line(1,0){10}}
\multiput(0,0)(10,0){2}{\line(0,1){10}}
\put(6.35,8.7){\circle*{0.4}}
\put(5.85,7.15){\circle*{0.4}}
\put(4.35,5.15){\circle*{0.4}}
\put(3.35,4.15){\circle*{0.4}}
\put(1.85,2.65){\circle*{0.4}}
\linethickness{0.7pt}
\put(2.5,3){\line(1,0){0.5}}\put(3,3){\line(0,1){1.5}}
\put(3,4.5){\line(1,0){1}}\put(4,4.5){\line(0,1){1}}
\put(4,5.5){\line(1,0){1.5}}\put(5.5,5.5){\line(0,1){2}}
\put(5.5,7.5){\line(1,0){0.5}}\put(6,7.5){\line(0,1){0.5}}
\bezier{5}(5.95,8)(5.95,8.5)(5.95,9)
\bezier{5}(1.5,3)(2,3)(2.5,3)
\end{picture}}
\hfill
\parbox[t]{11.75cm}
{Obviously, $\pi$ contains a decreasing subsequence of length $k$ if the 
diagram of $\pi$ has at least $k-1$ corners.\\ 
On the other hand, if there are at most $k-2$ corners in $D(\pi)$, we have $\des(\pi)\le k-2$ and hence $\pi\in{\cal 
S}_n(k\cdots1)$. (Note that each corner of $D(\pi)$ corresponds to a descent in $\pi$.)}
\end{minipage}
\vspace*{0cm}

{\bf b)} If the diagram $(n+1-k,n-k,\ldots,1)$ fits in $D(\pi)$ then 
we have $\pi_i\ge n+3-k-i$ for all $i$. Hence any increasing subsequence of $\pi$ is 
of length at most $k-1$: if $\pi_i=n+3-k-i$ then $i-1$ elements from $i+k-3$ 
many ones in $\{n+4-k-i,n+5-k-i,\ldots,n\}$ appear in $\pi_1\cdots\pi_{i-1}$.\\ 
Conversely, let $i$ be the smallest integer with $i+\pi_i<n+3-k$. Furthermore, 
choose $j$ such that $\pi_j=n+3-k-i$ (by definition of $i$, we have $j>i$), and let 
$\pi_{i_1}<\ldots<\pi_{i_{k-2}}$ be the elements of 
$\{n+4-k-i,n+5-k-i,\ldots,n\}$ which are not equal to $\pi_1,\ldots,\pi_{i-1}$. 
Note that $j<i_1<\ldots<i_{k-2}$ since $\pi$ is $132$-avoiding. Thus
$\pi_i\pi_j\pi_{i_1}\cdots\pi_{i_{k-2}}$ is an increasing sequence.
\vspace*{0.5ex}

{\bf c)}\hspace*{0.5cm} 
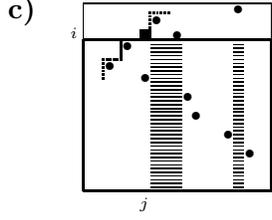
\begin{minipage}[t]{15cm}
\raisebox{-2.3cm}
{\unitlength=0.25cm
\begin{picture}(10,10)
\linethickness{0.2pt}
\multiput(0,0)(0,10){2}{\line(1,0){10}}
\multiput(0,0)(10,0){2}{\line(0,1){10}}
\put(5,8.3){\circle*{0.4}}\put(3.9,9.1){\circle*{0.4}}
\put(8.25,9.7){\circle*{0.4}}\put(2.35,7.7){\circle*{0.4}}
\put(1.4,6.65){\circle*{0.4}}\put(3.3,6){\circle*{0.4}}
\put(5.55,5){\circle*{0.4}}\put(6,4){\circle*{0.4}}
\put(7.7,3){\circle*{0.4}}\put(8.85,2){\circle*{0.4}}
\multiput(3.6,0)(0,0.2){40}{\line(1,0){1.6}}
\multiput(8,0)(0,0.2){40}{\line(1,0){0.5}}
\linethickness{0.7pt}
\put(1.95,8.05){\line(1,0){1.1}}\put(2,7){\line(0,1){1}}
\bezier{5}(1,7)(1.5,7)(1.95,7)\bezier{5}(1,7)(1,6.5)(1,6)
\bezier{5}(3.5,9.5)(4,9.5)(4.5,9.5)\bezier{5}(3.5,8.5)(3.5,9)(3.5,9.5)
\put(3,8){\rule{0.15cm}{0.15cm}}
\put(-0.7,8.4){\makebox(0,0)[lc]{\tiny$i$}}
\put(3,-0.7){\makebox(0,0)[lc]{\tiny$j$}}
\put(0,8.05){\line(1,0){10}}\put(0,0){\line(1,0){10}}
\put(0,0){\line(0,1){8.05}}\put(10,0){\line(0,1){8.05}}
\end{picture}}
\hfill
\parbox[t]{11.75cm}{Let $(i,j)$ be the top corner of $D(\pi)$ for which 
$i+j<n+3-k$. By removing the rows $1,\ldots,i$ and the columns 
$\pi_1,\ldots,\pi_i$, we obtain the diagram of a permutation $\sigma\in{\cal 
S}_{n-i}(132)$ whose letters are in the same relative order as 
$\pi_{i+1}\cdots\pi_n$ where $\sigma_1=\pi_{i+1}\le j<\pi_i$. As discussed in 
the proof of part b), the element $\sigma_1$ is the first one of an increasing}
\end{minipage}
\vspace*{-0.15cm}

sequence of length $(n-i)-\sigma_1+1$ in $\sigma$. (Since $\sigma_1\le 
j<n+3-k-i$ the index $l=1$ is the smallest one with $l+\sigma_l<(n-i)+3-k$.) Clearly, the 
first $j+1-\sigma_1$ terms are restricted by $j$. Thus there is an increasing sequence of 
length $n-(i+j)>k-3$ in $\sigma$ whose all elements are larger than $\pi_i$. 
Note that the elements $j+1,j+2,\ldots,\pi_i-1$ appear in 
$\pi_1\cdots\pi_{i-1}$.\\
To prove the converse, suppose that every corner of $D(\pi)$ satisfies the 
condition given above. Then we have $\pi_i+i>n+3-k$ for all $i\in{\sf D}(\pi)$. 
Hence for each descent $i$ of $\pi$ there exist at most $k-3$ elements $\pi_j$ 
with $j>i$ and $\pi_j>\pi_i$. Since $\pi$ is $132$-avoiding these elements form an increasing 
sequence. Thus there is no pattern $2134\cdots k$ in 
$\pi$.
\end{bew}

\begin{beme} \label{length of de/increasing subsequence}
\brm
\begin{enum}
\item[]
\item From the statement of a), it follows that the maximum length of a decreasing subsequence of 
$\pi\in{\cal S}_n(132)$ is equal to the number of corners of $D(\pi)$ plus one, 
or in terms of permutation statistics, $\des(\pi)+1$. It is well-known that 
the length of the longest decreasing sequence can easily be obtained for any 
permutation in ${\cal S}_n$ via the Robinson-Schensted correspondence: it is just 
the number of rows of one of the tableaux $P$ and $Q$ corresponding to 
$\pi\in{\cal S}_n$.\\
It is clear from the construction that every left-to-right-minimum of $\pi$ 
appears in the first column of $P$. Any other entry of this column must be 
the largest element of a $132$-pattern in $\pi$. 
Hence in case of $132$-avoiding permutations the elements of the first 
column of $P$ are precisely the left-to-right-minima. As observed in Section 2, these minima are just the 
descent bottoms and the first letter of $\pi$. Thus for $\pi\in{\cal 
S}_n(132)$ the tableau $P$ has $\des(\pi)+1$ rows. 
\item By part b), the length of a longest increasing subsequence of 
$\pi\in{\cal S}_n(132)$ equals the maximum value of $n+1-i-\lambda_i$ with 
$1\le i\le n-1$ where $\lambda_i\ge0$ is the length of the $i$th row of $D(\pi)$. This also follows from Remark \ref{remarks to section 
2}b) according to which $\lambda_i$ is equal to the $i$th component of the code 
of $\pi$. Thus $n-i-\lambda_i$ counts the number of elements on the right of $\pi_i$ 
which are larger than $\pi_i$. (Since $\pi$ contains no pattern $132$ these 
elements appear in increasing order.)
\end{enum}
\erm
\end{beme}

\begin{kor}
$|{\cal S}_n(132,k(k-1)\cdots1)|=\frac{1}{n}\sum_{i=1}^{k-1}{n\choose 
i}{n\choose i-1}$ for all $n$ and $k\ge3$.
\end{kor}

\begin{bew}
As mentioned in \ref{remarks to section 2}c) the number of partitions 
in ${\cal Y}_n$ whose diagram has exactly $i$ corners is equal to the Narayana 
number $N(n,i+1)=\frac{1}{n}{n\choose i}{n\choose i+1}$. Thus there are 
$\sum_{i=0}^{k-2} N(n,i+1)$ diagrams with at most $k-2$ corners.
\end{bew}
\vspace*{1ex}

The following result also follows from a special case of \cite[Th. 2.6]{mansour-vainshtein3}.

\begin{kor}
$|{\cal S}_n(132,12\cdots k)|=|{\cal S}_n(132,213\cdots k)|$ for all $n$ and 
$k\ge3$.
\end{kor}

\begin{bew}
There is a simple bijection between the restricted diagrams which contain 
$(n+1-k,n-k,\ldots,1)$ and those ones whose all corners satisfy the 
condition $i+j\ge n+3-k$. (Note that the empty diagram associating with the 
identity in ${\cal S}_n$ belongs to the latter ones.) For each corner $(i,j)$ of the diagram 
$(n+1-k,n-k,\ldots,1)$ we have $i+j=n+2-k$. Thus every diagram containing $(n+1-k,n-k,\ldots,1)$ 
is uniquely determined by its corners outside this shape (which are precisely 
the corners with $i+j\ge n+3-k$). Given such a diagram $D$, the corresponding 
diagram $D'$ is defined to be this one whose corners are the corners of $D$ 
which are not contained in $(n+1-k,n-k,\ldots,1)$. Conversely, for any diagram 
$D'$ whose all corners satisfy $i+j\ge n+3-k$ we construct the corresponding 
diagram $D$ as the union of $D'$ and $(n+1-k,n-k,\ldots,1)$.
\end{bew}
\vspace*{1ex}

While the above conditions can be checked without effort, the characterization 
of the avoidance of the patterns considered now is more technical.\\[2ex]
Given a permutation $\pi\in{\cal S}_n(132)$, let 
$\lambda_1,\ldots,\lambda_l$ be the positive parts of the partition 
with the diagram $D(\pi)$. Let $a_i:=n-(i+\lambda_i)$ for $i=1,\ldots,l$ and 
$b_i:=n-(i+\lambda'_i)$ for $i=1,\ldots,\lambda_1$ where $\lambda'$ denotes the 
conjugate of $\lambda$. Furthermore, for $i=1,\ldots,l$ let $h_i$ be the length of the longest 
increasing sequence in $b_{\lambda_i}b_{\lambda_i-1}\cdots b_1$ whose first 
element is $b_{\lambda_i}$. We call 
the number $h_i$ the {\it height of $a_i$}.\\
For example, the permutation $\pi=8\:9\:5\:4\:6\:7\:2\:3\:10\:1\in{\cal 
S}_{10}(132)$ generates the diagram of $\lambda=(7,7,4,3,3,3,1,1,1)$:
\vspace*{-1ex}

\begin{center}
\unitlength=0.3cm
\begin{picture}(40,11)
\linethickness{0.2pt}
\put(1,0){\line(0,1){2}}\put(11,0){\line(0,1){10}}
\put(1,0){\line(1,0){10}}\put(8,10){\line(1,0){3}}
\put(2,10){\line(0,-1){8}}\put(3,10){\line(0,-1){6}}
\put(4,10){\line(0,-1){3}}\multiput(5,10)(1,0){3}{\line(0,-1){2}}
\put(1,9){\line(1,0){7}}\put(1,8){\line(1,0){4}}
\multiput(1,7)(0,-1){3}{\line(1,0){3}}\multiput(1,4)(0,-1){3}{\line(1,0){1}}
\linethickness{0.7pt}
\put(1,1){\line(0,1){9}}\put(2,1){\line(0,1){3}}
\put(1,10){\line(1,0){7}}\put(1,1){\line(1,0){1}}
\put(2,4){\line(1,0){2}}\put(4,4){\line(0,1){3}}
\put(4,7){\line(1,0){1}}\put(5,7){\line(0,1){1}}
\put(5,8){\line(1,0){3}}\put(8,8){\line(0,1){2}}
\put(8.5,9.5){\circle*{0.4}}\put(9.5,8.5){\circle*{0.4}}
\put(5.5,7.5){\circle*{0.4}}\put(4.5,6.5){\circle*{0.4}}
\put(6.5,5.5){\circle*{0.4}}\put(7.5,4.5){\circle*{0.4}}
\put(2.5,3.5){\circle*{0.4}}\put(3.5,2.5){\circle*{0.4}}
\put(10.5,1.5){\circle*{0.4}}\put(1.5,0.5){\circle*{0.4}}
\put(0.5,9.5){\makebox(0,0)[cc]{\sf\tiny8}}
\put(0.5,8.5){\makebox(0,0)[cc]{\sf\tiny9}}
\put(0.5,7.5){\makebox(0,0)[cc]{\sf\tiny7}}
\put(0.5,6.5){\makebox(0,0)[cc]{\sf\tiny7}}
\put(0.5,5.5){\makebox(0,0)[cc]{\sf\tiny8}}
\put(0.5,4.5){\makebox(0,0)[cc]{\sf\tiny9}}
\put(0.5,3.5){\makebox(0,0)[cc]{\sf\tiny8}}
\put(0.5,2.5){\makebox(0,0)[cc]{\sf\tiny9}}
\put(0.32,1.5){\makebox(0,0)[cc]{\sf\tiny10}}
\put(1.5,10.5){\makebox(0,0)[cc]{\sf\tiny10}}
\put(2.5,10.5){\makebox(0,0)[cc]{\sf\tiny8}}
\put(3.5,10.5){\makebox(0,0)[cc]{\sf\tiny9}}
\put(4.5,10.5){\makebox(0,0)[cc]{\sf\tiny7}}
\put(5.5,10.5){\makebox(0,0)[cc]{\sf\tiny7}}
\put(6.5,10.5){\makebox(0,0)[cc]{\sf\tiny8}}
\put(7.5,10.5){\makebox(0,0)[cc]{\sf\tiny9}}
\put(13,1.5){\makebox(0,0)[lb]{\footnotesize The numbers $i+\lambda_i$ and 
$i+\lambda'_i$ are given on the left side}}
\put(13,0){\makebox(0,0)[lb]{\footnotesize and on the top of the diagram, 
respectively.}}
\end{picture}
\end{center}
So we obtain $a(\pi)=(2,1,3,3,2,1,2,1,0)$ and $h(\pi)=(3,3,1,2,2,2,1,1,1)$. 
(Note that we have $b(\pi)=(0,2,1,3,3,2,1)$.)

\begin{satz} \label{shifted patterns}
Let $\pi\in{\cal S}_n(132)$ be a permutation, and let $a(\pi),h(\pi)$ be as above. Then 
$\pi$ avoids the pattern $s(s+1)\cdots k12\cdots(s-1)$, where $2\le s\le k$ and 
$k\ge3$, if and only if the longest decreasing subsequence of $a(\pi)$ which ends in 
an element of height $\ge s-1$ is of length at most $k-s$.
\end{satz}

\begin{bew}
Let $i_1$ be an integer such that $a_{i_1}\ge a_i$ for all $i<i_1$. Hence a 
decreasing sequence whose first element is $a_{i_1}$ can not be extended to the 
left. Then $\lambda_{i_1}<\lambda_{i_1-1}$ where $\lambda_0:=n$. (Note that the 
condition $\lambda_i<\lambda_{i-1}$ is equivalent to $a_{i-1}\le a_i$.) As 
shown in Section 2, the element $\pi_{i_1}$ is a left-to-right minimum of 
$\pi$. Thus any increasing subsequence in $\pi$ which starts with $\pi_{i_1}$ 
is left maximal. In particular, we have $\pi_{i_1}=\lambda_{i_1}+1$.\\
Now let $a_j$ be an element of $a(\pi)$ with $a_j<a_{i_1}$, $j>i_1$, and 
$a_{i_1+1},\ldots,a_{j-1}>a_j$. Since $j-i>\lambda_i-\lambda_j$ for 
$i=i_1,i_1+1,\ldots,j-1$, all the elements 
$\lambda_j+1,\lambda_j+2,\ldots,\lambda_{i_1}$ occur in 
$\pi_{i_1+1}\pi_{i_1+2}\cdots\pi_{j-1}$. Hence $\pi_{i_1}<\pi_j$.\\ 
Consequently, if $a_{i_1}>a_{i_2}>\ldots>a_{i_r}$ (with $i_1<i_2<\ldots<i_r$) is a sequence such that there is no integer 
$i$ with $i_l<i<i_{l+1}$ and $a_{i_l}>a_i\ge a_{i_{l+1}}$ for any $l$, then 
$\pi_{i_1}<\pi_{i_2}<\ldots<\pi_{i_r}$ is an increasing subsequence of $\pi$ 
which is maximal with respect to the property that $\pi_{i_1}$ and $\pi_{i_r}$ 
are its first and last elements, respectively. (Since $\pi$ avoids the pattern 
$132$, the relations $\pi_{i_1}<\pi_{i_l}$ for $l=2,\ldots,r$ imply that 
$\pi_{i_1},\ldots,\pi_{i_r}$ are increasing.)\\
It is clear from the definition that $h_i$ is the maximal length of an 
increasing sequence of dots southwest of the dot representing $\pi_i$ which 
begins at the top dot southwest of $(i,\pi_i)$. Thus, if $a_{i_r}$ is an 
element of height $\ge s-1$ then there exist at least $s-1$ integers 
$i_r<j_1<j_2<\ldots<j_{s-1}$ with $\pi_{j_1}<\ldots<\pi_{j_{s-1}}<\pi_{i_r}$. 
Since $\pi$ is $132$-avoiding, we even have 
$\pi_{j_1}<\ldots<\pi_{j_{s-1}}<\pi_{i_1}$. Choosing $i_1$ and $i_r$ minimal 
and maximal, respectively, proves the assertion.
\end{bew}
\vspace*{1ex}

For any permutation $\pi\in{\cal S}_n(132)$ denote by $l_s(\pi)$ the largest 
integer $l$ such that $\pi$ contains the shifted pattern $s(s+1)\cdots l12\cdots (s-1)$ 
where $s\ge2$. By the theorem, $l_s(\pi)$ is equal to 
$s-1$ plus the maximum length of a decreasing sequence in $a(\pi)$ whose 
smallest element is of height at least $s-1$.\\
It is clear that the sequence $L(\pi):=(l_2(\pi)-1,l_3(\pi)-2,\ldots)$ is a partition, 
that is, $l_s(\pi)+1\ge l_{s+1}(\pi)$ for all $s$. Since $\pi$ avoids $s(s+1)\cdots k12\cdots (s-1)$ if and only if no pattern 
$(k+2-s)$\\$(k+3-s)\cdots k12\cdots (k+1-s)$ occurs in $\pi^{-1}$, the partition 
$L(\pi^{-1})$ is the conjugate of $L(\pi)$. (Obviously, for any permutation $\pi\in{\cal S}_n$ the diagram of the inverse 
$\pi^{-1}$ is just the transpose of $D(\pi)$. Hence 
the set ${\cal S}_n(132)$ is closed under inversion.) 

\begin{bem} \label{shifted patterns in terms of Dyck path}
\brm
Using the relation between the diagram $D(\pi)$ and the Dyck path $\Psi_K(\pi)$, it is easy to see that the 
number $a_i$ is just the height at which the $i$th down-step of $\Psi_K(\pi)$ 
ends. (We only consider the down-steps before the last up-step.) The numbers $b_i$ 
needed for the construction of $h(\pi)$ are (in 
reverse order) the starting heights of the up-steps after the first down-step. 
Thus the height $h_i$ is precisely the number of (not necessarily consecutive) up-steps of increasing 
starting heights following the $i$th down-step in $\Psi_K(\pi)$. Hence, in case of $s=2$ the theorem yields the second part of \cite[Lemma $\Phi$]{krattenthaler}.
\erm
\end{bem}

We shall prove now that the number of permutations in ${\cal S}_n$ which avoid 
both $132$ and the pattern $s(s+1)\cdots k12\cdots (s-1)$ with $k\ge 3$ and 
$s\in\{1,\ldots,k\}$ does {\it not} depend on $s$.

\begin{prop} \label{case s=2}
Let $\pi\in{\cal S}_n(132)$, and let $l$ be the maximum length of an increasing 
subsequence of $\pi$. Then $\pi$ corresponds in a one-to-one fashion to a 
permutation $\sigma\in{\cal S}_n(132)$ with $l_2(\sigma)=l$.
\end{prop}

\begin{bew}
Let $\lambda$ and $\mu$ be the partitions whose diagrams equal $D(\pi)$ and 
$D(\sigma)$, respectively. Given $\lambda=(\lambda_1,\ldots,\lambda_{n-1})\in{\cal 
Y}_n$, we define the sequence $\hat{\mu}$ by
\bdpm
\hat{\mu}_i:=\left\{\ba{ccl}\lambda_i+1&&\mbox{if }\lambda_i+i<n\\
0&&\mbox{if }\lambda_i+i=n\ea\right.
\edpm
for $i=1,\ldots,n-1$, and obtain the partition $\mu$ by sorting $\hat{\mu}$. (Delete all parts 
$\hat{\mu}_i=0$ with nonzero $\hat{\mu}_{i+1}$ and add the corresponding number 
of zeros at the end of the sequence.) It is obvious that $\mu\in{\cal Y}_n$, 
and it is easy to see that the map $\lambda\mapsto\mu$ is injective, and 
thus a bijection on ${\cal Y}_n$. To recover $\lambda$ from $\mu$, first set 
$\hat{\lambda}_i:=\mu_i-1$ for all positive $\mu_i$. 
Then for any $\mu_i=0$ let $j$ be the largest integer for which 
$\hat{\lambda}_j+j\ge n-1$, and define $\hat{\lambda}$ to be the sequence obtained 
by inserting $n-1-j$ between $\hat{\lambda}_j$ and $\hat{\lambda}_{j+1}$. If 
there is no such an integer $j$ then prepend $n-1$ to $\hat{\lambda}$. The 
partition resulting from this procedure is defined to be $\lambda$.\\ 
Consider now the sequence $\bar{a}(\pi)=(n-i-\lambda_i)_{i=1,\ldots,n-1}$. (For the statement of Theorem \ref{shifted patterns} it 
suffices to consider the reduced sequence $a(\pi)$ which is obtained by 
omitting the final terms $\bar{a}_i=n-i$. By Remark \ref{length of de/increasing 
subsequence}b), we have $l=\max\bar{a}_i+1$. Let $j_1$ be an integer satisfying $\bar{a}_{j_1}=l-1$. Using the 
definition of $\bar{a}(\pi)$ in terms of the Dyck path $\Psi_K(\pi)$, it is 
obvious that there exist some integers $j_1<j_2<\ldots<j_{l-1}\le n-1$ with $\bar{a}_{j_i}=l-i$ for 
$i=1,\ldots,l-1$. By Remark \ref{shifted patterns in terms of Dyck 
path}, the number $\bar{a}_i(\pi)$ is the height of the ending point of the 
$i$th down-step of $\Psi_K(\pi)$. (Note that the length $l$ is equal to the maximum height of a 
peak of $\Psi_K(\pi)$.) By the construction, the elements of the sequence 
$a(\sigma)$ correspond to the nonzeros of $\bar{a}(\pi)$. In particular, we 
have $a_i(\sigma)>a_j(\sigma)$ for any $i<j$ if and only if the $i$th positive 
element of $\bar{a}(\pi)$ is larger than the $j$th positive one. Thus the 
elements corresponding to $\bar{a}_{j_1},\bar{a}_{j_2},\ldots,\bar{a}_{j_{l-1}}$ form a 
decreasing sequence of maximal length in $a(\sigma)$. (For $l=1$ we have 
$\lambda=(n-1,n-2,\ldots,1)$, i.e., $\pi=n(n-1)\cdots1$ and hence 
$\mu=\emptyset$, i.e., $\sigma={\rm id}\in{\cal S}_n$.) Consequently, 
$l_2(\sigma)=(l-1)+1$. (Note that any element 
$a_i(\sigma)$ is of height at least 1.) 
\end{bew}
\vspace*{1ex}

The following result can be derived from the corresponding generating functions 
which were given for the first time by Chow and West (\cite[Th. 3.1]{chow-west}). 
Several different analytical proofs appeared recently in 
\cite[Th. 2, Th. 6]{krattenthaler} and \cite[Th. 3.1]{mansour-vainshtein2}.

\begin{kor} \label{12...k=23...k1}
$|{\cal S}_n(132,12\cdots k)|=|{\cal S}_n(132,23\cdots k1)|=|{\cal S}_n(132,k12\cdots 
(k-1))|$ for all $n$ and $k\ge3$.
\end{kor}

\begin{bew}
The first identity is an immediate consequence of the preceding proposition. For 
the second one use that $\pi$ avoids $23\cdots k1$ if and only if $\pi^{-1}$ 
contains no pattern $k12\cdots(k-1)$.
\end{bew}

\begin{kor} \label{l smaller s}
For any $s\ge2$ there are as many partitions $\lambda\in{\cal Y}_n$ for which 
$i+\lambda_i\ge n+2-s$ for all $i$ as such ones for which the sequence $a(\pi)$ 
for the corresponding permutation contains no element of height at least $s-1$.
\end{kor}

\begin{bew}
If $\lambda\in{\cal Y}_n$ satisfies the condition $i+\lambda_i\ge n+2-s$ for 
all $i$ then its diagram contains $(n+1-s,n-s,\ldots,1)$. By Theorem 
\ref{Diagrammcharakterisierung}b), this is equivalent to being $12\cdots 
s$-avoiding for the corresponding permutation. On the other hand, it follows from Theorem \ref{shifted patterns} 
(where $k=s$) that any permutation $\pi\in{\cal S}_n(132)$ avoids the pattern $s12\cdots(s-1)$ if and only if every element of 
$a(\pi)$ is of height at most $s-2$. Thus Corollary \ref{12...k=23...k1} proves the 
assertion. (Let $\tau$ be the partition associating with $\pi$. Then 
$\lambda\mapsto\mu\mapsto\tau:=\mu'$ yields the desired bijection where 
$\lambda\mapsto\mu$ is the map appearing in the proof of \ref{case s=2}, and 
$\mu'$ denotes the conjugate of $\mu$.)
\end{bew}
\vspace*{1ex}

Now we can generalize the result from Proposition \ref{case s=2} for any $s\ge2$. 

\begin{prop} 
Let $\pi\in{\cal S}_n(132)$, and let $l\ge s-1$ be the maximum length of an increasing 
subsequence of $\pi$. Then $\pi$ corresponds in a one-to-one fashion to a 
permutation $\sigma\in{\cal S}_n(132)$ with $l_s(\sigma)=l$.
\end{prop}

\begin{bew}
For $l\le s-1$ we have $i+\lambda_i\ge n+2-s$ for all parts $\lambda_i$ of the 
partition corresponding to $\pi$. (Note that $l$ 
is the maximum difference $n+1-(i+\lambda_i)$.) In this case, by Corollary \ref{l smaller 
s}, $\lambda$ corresponds to a partition whose associated permutation 
satisfies $l_s=s-1$. Thus we may assume that $l\ge s$.\\  
The argument is completely parallel to that used to prove Proposition \ref{case 
s=2} but the construction of the bijection is more complicated. We preserve the notation of 
the proof of \ref{case s=2}.\\ 
Given $\lambda=(\lambda_1,\ldots,\lambda_{n-1})\in{\cal 
Y}_n$, we define the sequence $\hat{\mu}$ by
\bdpm
\hat{\mu}_i:=\left\{\ba{ccl}\lambda_i+s-1&&\mbox{if }\lambda_i+i<n+2-s\\
0&&\mbox{if }\lambda_i+i=n+2-s\\\vdots&&\qquad\quad\vdots\\s-2&&\mbox{if }\lambda_i+i=n\\\ea\right.
\edpm
for $i=1,\ldots,n+1-s$, and $\hat{\mu}_i:=\lambda_i$ otherwise. If there exists 
any integer $i$ with $\hat{\mu}_i<s-1$ and $\hat{\mu}_{i+1}\ge s-1$ then sort 
$\hat{\mu}$ as follows: If $\hat{\mu}_i>0$ then increase $\hat{\mu}_{i+1}$ by 
$1$. Interchange $\hat{\mu}_i$ and $\hat{\mu}_{i+1}$. After this procedure we 
have $\hat{\mu}=\tau_1\tau_2$ where $\tau_1$ is a partition with parts at least $s-1$, and $\tau_2$ is a subsequence 
containing elements of $\{0,\ldots,s-2\}$ only. (Note that 
both $\tau_1$ and $\tau_2$ are nonempty for $s\ge3$. Since $l\ge s$ there is an 
integer $i$ satisfying $i+\lambda_i<n+2-s$ and hence $\hat{\mu}_i\ge s-1$. On 
the other hand, $\hat{\mu}_{n+2-s},\ldots,\hat{\mu}_{n-1}\le s-2$ by 
definition.) It is clear from the construction that $i+\hat{\mu}_i\le n$, and 
it is easy to see that the map $\lambda\mapsto\hat{\mu}$ is injective. To 
recover $\lambda$ from $\hat{\mu}$, first define a sequence $\hat{\lambda}$ by 
$\hat{\lambda}_i:=\hat{\mu}_i+1-s$ for all indices $i$ with $\hat{\mu}_i\ge s-1$ and 
$\hat{\lambda}_i:=\hat{\mu}_i-i+n+2-s$ otherwise where $i=1,\ldots,n+1-s$. If 
there exists any integer with $\hat{\lambda}_i<\hat{\lambda}_{i+1}$ then 
decrease $\hat{\mu}_i$ by 1 if $\hat{\mu}_{i+1}$ is positive, interchange 
$\hat{\mu}_i$ and $\hat{\mu}_{i+1}$, and determine $\hat{\lambda}$ again for the 
resulting sequence $\hat{\mu}$. (This procedure will be done while 
$\hat{\lambda}$ is not a partition.) Adding 
$\hat{\mu}_{n+2-s},\ldots,\hat{\mu}_{n-1}$ at $\hat{\lambda}$ completes the 
partition $\lambda$.\\
If $\tau_2$, and hence $\hat{\mu}$, is a partition then set $\mu:=\hat{\mu}$. The 
partitions $\lambda\in{\cal Y}_n$ for which this case occurs are characterized by the 
following conditions:
\begin{enum2}
\item $\lambda_{n+2-s}=\ldots=\lambda_{n-1}=0$. 
\item If $i+\lambda_i\ge n+2-s$ for any $i\le n-s$ then 
$\lambda_{i+1}<\lambda_i$. 
\end{enum2}
\vspace*{-0.35cm}

Suppose that there exists an integer $i\ge n+2-s$ with $\hat{\mu}_i=\lambda_i>0$. 
By Remark \ref{shifted patterns in terms of Dyck path}, the $i$th down-step of 
$\Psi_K(\pi)$ ends at a point of height $n-(i+\lambda_i)\le 
n-(n+2-s+\lambda_i)=s-2-\lambda_i\le s-3$. Since $l\ge s$ there is an 
integer $j$ with $j+\lambda_j<n+2-s$. (Obviously, $j<n+2-s\le i$.) The ending point 
of the corresponding down-step is of height at least $s-1$. Thus there must 
exist an integer $k$ with $j<k<i$ such that $k+\lambda_k=n+2-s$, that is, 
$\hat{\mu}_k=0$. The second condition is evident. (In case $s=2$ the above 
conditions are satisfied for all partitions in ${\cal Y}_n$.)\\
If $\tau_2$ is not ordered then we obtain the partition $\mu$ by increasing all 
elements of $\hat{\mu}$ which are smaller than $s-1$ and smaller than an 
element to their right.\\
Consider the sequences $\bar{a}(\pi)$ and $a(\sigma)$ again. By the 
construction, each element of $\bar{a}(\pi)$ which is at least $s-1$ corresponds to a term $a_i$ of $a(\sigma)$ 
whose height is at least $s-1$ or for which there exists an element $a_j$ with 
$i<j$ and $h_j(\sigma)\ge s-1$. (Then necessarily $a_i\ge a_j$.) In 
particular, we have $a_{i_1}(\sigma)>a_{i_2}(\sigma)$ for any $i_1<i_2\le i_3$ 
where $h_{i_3}(\sigma)\ge s-1$ if and only if 
the corresponding elements in $\bar{a}(\pi)$ are in the same order. Analogously 
to the proof of Proposition \ref{case s=2} this yields the assertion.
\end{bew}
\vspace*{1ex}

An analytical proof of the following result was given in 
\cite[Th. 2.4]{mansour-vainshtein3}.

\begin{kor}
$|{\cal S}_n(132,12\cdots k)|=|{\cal S}_n(132,s(s+1)\cdots k12\cdots 
(s-1))|$ for all $n$ and $k\ge 3$ and $2\le s\le k$.
\end{kor}

\begin{bew}
This follows immediately from the previous proposition. 
\end{bew}
\vspace*{0.5cm}


\setcounter{section}{5}\setcounter{satz}{0}

\centerline{\large{\bf 5}\hspace*{0.25cm}
{\sc A final note}}
\vspace*{0.5cm}

In Section 2 we have shown that the diagram of a permutation indicates the 
existence of some subsequences of type $132$. But in positive case, we even obtain the exact number of 
occurrences.\\[2ex]
In \cite{fulton}, Fulton defined the following rank function on the essential 
set. Given a corner $(i,j)$ of the diagram $D(\pi)$, i.e. $(i,j)\in{\cal 
E}(\pi)$, its {\it rank} is defined by the number of dots northwest of $(i,j)$ 
and denoted by $\rho((i,j))$.\\
It is clear from the construction that the number of dots in the northwest is the same 
for all diagram squares which are connected. (This yields the fundamental 
property of the ranked essential set of a permutation $\pi$, that it uniquely determines 
$\pi$.) So, we can extend the rank function on $D(\pi)$. The information about 
the number of sequences of type $132$ containing in a permutation 
is encoded by the ranks of its diagram squares.

\begin{satz}
Let $\pi\in{\cal S}_n$ be a permutation, and let $D(\pi)$ be its diagram. Then 
the number of occurrences of the pattern $132$ in $\pi$ is equal to
\bdpm
\sum_{(i,j)\in D(\pi)} \rho((i,j)).
\edpm
\end{satz}
\vspace*{0.3cm}

\begin{bew}
Extending the arguments of Theorem \ref{132-avoiding diagram}, it is easy 
to see that each square $(i,j)$ of $D(\pi)$ corresponds to exactly 
$\rho((i,j))$ subsequences of type $132$ in $\pi$, namely the sequences 
$k\:\pi_i\:j$ where $k$ ranges over all column indices of dots northwest of $(i,j)$:
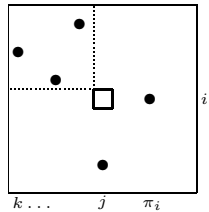
\begin{center}
\unitlength=0.25cm
\begin{picture}(10,10)
\linethickness{0.2pt}
\multiput(0,0)(0,10){2}{\line(1,0){10}}
\multiput(0,0)(10,0){2}{\line(0,1){10}}
\put(3.75,9){\circle*{0.5}}
\put(0.5,7.5){\circle*{0.5}}
\put(2.5,6){\circle*{0.5}}
\put(5,1.5){\circle*{0.5}}
\put(7.5,5){\circle*{0.5}}
\put(10.2,5){\makebox(0,0)[lc]{\tiny$i$}}
\put(5,-0.5){\makebox(0,0)[cc]{\tiny$j$}}
\put(0.2,-0.5){\makebox(0,0)[lc]{\tiny$k\ldots$}}
\put(7.1,-0.65){\makebox(0,0)[lc]{\tiny$\pi_i$}}
\linethickness{0.3pt}
\bezier{20}(4.5,5.5)(4.5,7.75)(4.5,9.9)
\bezier{20}(4.5,5.5)(2.25,5.5)(0.1,5.5)
\linethickness{0.7pt}
\put(4.5,5.5){\line(1,0){1}}\put(4.5,4.5){\line(1,0){1}}
\put(4.5,4.5){\line(0,1){1}}\put(5.5,4.5){\line(0,1){1}}
\end{picture}
\end{center}
\vspace*{-0.5cm}
\end{bew}

\begin{bem}
\brm
As mentioned above, we have $|D(\pi)|=\inv(\pi)$ for all $\pi\in{\cal S}_n$. 
Hence the non-weighted sum $\sum_{(i,j)\in D(\pi)} 1$ counts the number of occurrences of the 
pattern $21$ in $\pi$.
\erm
\end{bem}

\begin{bsp}
\brm
The ranked diagram of $\pi=4\:2\:8\:3\:6\:9\:7\:5\:1\:10\in{\cal S}_{10}$ is
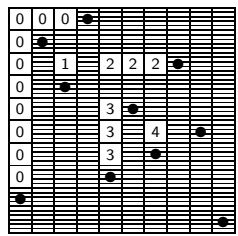
\begin{center}
\unitlength=0.3cm
\begin{picture}(10,10)
\linethickness{0.3pt}
\multiput(0,0)(0,1){11}{\line(1,0){10}}
\multiput(0,0)(1,0){11}{\line(0,1){10}}
\put(3.5,9.5){\circle*{0.5}}
\put(1.5,8.5){\circle*{0.5}}
\put(7.5,7.5){\circle*{0.5}}
\put(2.5,6.5){\circle*{0.5}}
\put(5.5,5.5){\circle*{0.5}}
\put(8.5,4.5){\circle*{0.5}}
\put(6.5,3.5){\circle*{0.5}}
\put(4.5,2.5){\circle*{0.5}}
\put(0.5,1.5){\circle*{0.5}}
\put(9.5,0.5){\circle*{0.5}}
\linethickness{0.2pt}
\multiput(3,9)(0,0.2){5}{\line(1,0){7}}
\multiput(1,8)(0,0.2){5}{\line(1,0){9}}
\multiput(7,7)(0,0.2){5}{\line(1,0){3}}
\multiput(2,6)(0,0.2){5}{\line(1,0){8}}
\multiput(5,5)(0,0.2){5}{\line(1,0){5}}
\multiput(8,4)(0,0.2){5}{\line(1,0){2}}
\multiput(6,3)(0,0.2){5}{\line(1,0){4}}
\multiput(4,2)(0,0.2){5}{\line(1,0){6}}
\multiput(0,1)(0,0.2){5}{\line(1,0){10}}
\multiput(0,0)(0,0.2){5}{\line(1,0){10}}
\multiput(1,2)(0,0.2){5}{\line(1,0){3}}
\multiput(1,3)(0,0.2){5}{\line(1,0){3}}
\multiput(5,3)(0,0.2){5}{\line(1,0){1}}
\multiput(1,4)(0,0.2){5}{\line(1,0){3}}
\multiput(1,5)(0,0.2){5}{\line(1,0){3}}
\multiput(1,6)(0,0.2){5}{\line(1,0){1}}
\multiput(1,7)(0,0.2){5}{\line(1,0){1}}
\multiput(3,7)(0,0.2){5}{\line(1,0){1}}
\multiput(5,4)(0,0.2){5}{\line(1,0){1}}
\multiput(7,4)(0,0.2){5}{\line(1,0){1}}
\multiput(0.5,9.5)(1,0){3}{\makebox(0,0)[cc]{\sf\tiny0}}
\multiput(0.5,8.5)(0,-1){7}{\makebox(0,0)[cc]{\sf\tiny0}}
\multiput(4.5,7.5)(1,0){3}{\makebox(0,0)[cc]{\sf\tiny2}}
\multiput(4.5,5.5)(0,-1){3}{\makebox(0,0)[cc]{\sf\tiny3}}
\put(2.5,7.5){\makebox(0,0)[cc]{\sf\tiny1}}
\put(6.5,4.5){\makebox(0,0)[cc]{\sf\tiny4}}
\end{picture}
\end{center}
Thus $\pi$ contains $20$ subsequences of type $132$ and 18 inversions.
\erm
\end{bsp}

In \cite{mansour-vainshtein4}, Mansour and Vainshtein studied the generating 
function for the number of permutations on $n$ letters containing exactly $r\ge 
0$ occurrences of pattern $132$.
\vspace*{1cm}

{\bf Acknowledgement} I am grateful to Toufik Mansour for pointing me to some 
literature relating to restricted permutations.
\newpage


\centerline{\large\sc References}
\vspace*{0.5cm}

\begin{enumbib}

\bibitem{billey etal}
S. C. Billey, W. Jockusch and R. P. Stanley, 
Some Combinatorial Properties of Schubert Polynomials, 
{\it J. Alg. Comb.} {\bf 2} (1993), 345-374.

\bibitem{braenden etal}
P. Br\"and\'en, A. Claesson and E. Steingr\'{\i}msson, 
Catalan continued fractions and increa\-sing subsequences in permutations, 
{\it Discrete Math.} {\bf 258} (2002), 275-287. 

\bibitem{chow-west}
T. Chow and J. West, 
Forbidden subsequences and Chebyshev polynomials, 
{\it Discrete Math.} {\bf 204} (1999), 119-128. 

\bibitem{deutsch}
E. Deutsch, 
A bijection on Dyck paths and its consequences, 
{\it Discrete Math.} {\bf 179} (1998), 253-256. 

\bibitem{eriksson-linusson}
K. Eriksson and S. Linusson, 
Combinatorics of Fulton's essential set,
{\it Duke Math. J.} {\bf 85} (1996), 61-80.

\bibitem{feller}
W. Feller, 
An Introduction to Probability Theory and Its Applications, vol. I,
Wiley, New York, 1968.

\bibitem{fulmek}
M. Fulmek, 
Enumeration of permutations containing a prescribed number of occurrences of a pattern of length 3, 
preprint math.CO/0112092, 2001.

\bibitem{fulton}
W. Fulton, 
Flags, Schubert polynomials, degeneracy loci, and determinantal formulas,
{\it Duke Math. J.} {\bf 65} (1992), 381-420.

\bibitem{krattenthaler}
C. Krattenthaler, 
Permutations with restricted patterns and Dyck paths,
{\it Adv. Appl. Math.} {\bf 27} (2001), 510-530. 

\bibitem{macdonald}
I. G. Macdonald, 
Notes on Schubert Polynomials,
LaCIM, Universit\'{e} du Qu\'{e}bec \`{a} Montr\'{e}al, 1991.

\bibitem{mansour-vainshtein1}
T. Mansour and A. Vainshtein, 
Restricted permutations and Chebyshev polynomials, 
{\it S\'{e}min. Lothar. Comb.} {\bf 47} (2002), B47c.

\bibitem{mansour-vainshtein2}
T. Mansour and A. Vainshtein, 
Restricted permutations, continued fractions, and Chebyshev polynomials, 
{\it Electron. J. Comb.} {\bf 7} (2000), R17.

\bibitem{mansour-vainshtein3}
T. Mansour and A. Vainshtein, 
Restricted $132$-avoiding permutations, 
{\it Adv. Appl. Math.} {\bf 26} (2001), 258-269. 

\bibitem{mansour-vainshtein4}
T. Mansour and A. Vainshtein, 
Counting occurrences of $132$ in a permutation, 
{\it Adv. Appl. Math.} {\bf 28} (2002), 185-195. 

\bibitem{reifegerste}
A. Reifegerste, 
Differenzen in Permutationen: \"Uber den Zusammenhang von Permutationen, 
Polyominos und Motzkin-Pfaden,  
Ph.D. Thesis, University of Magdeburg, 2002.

\bibitem{reifegerste1}
A. Reifegerste, 
The excedances and descents of bi-increasing permutations, in preparation.

\bibitem{simion-schmidt}
R. Simion and F. W. Schmidt, 
Restricted Permutations, 
{\it Europ. J. Combinatorics} {\bf 6} (1985), 383-406. 

\bibitem{west}
J. West, 
Generating trees and the Catalan and Schr\"oder numbers,
{\it Discrete Math.} {\bf 146} (1995), 247-262. 
 
\end{enumbib}
 
\end{document}